\documentstyle[11pt]{article}
\def\P{{\bf P}}
\def\C{{\cal C}}
\def\K{{\cal K}}
\def\O{{\cal O}}

\begin{document}

\begin{center}
{\Large \bf Unobstructedness of deformations of holomorphic maps }\\
{\Large \bf onto Fano manifolds of Picard number 1}

\bigskip
 {\large \bf Jun-Muk Hwang}\footnote{This work was supported by the SRC Program ASARC
  funded by the Korea government.}

 \end{center}

\bigskip
\begin{abstract}
We show that deformations of a surjective morphism onto a Fano
manifold of Picard number 1 are unobstructed and  rigid modulo the
automorphisms of the target, if the variety of minimal rational
tangents of the Fano manifold is non-linear or finite. The
condition on the variety of minimal rational tangents holds for
practically all known examples of Fano manifolds of Picard number
1, except the projective space. When the variety of minimal
rational tangents is non-linear, the proof is based on an earlier
result of N. Mok and the author on the birationality of the
tangent map. When the varieties of minimal rational tangents of
the Fano manifold is finite, the key idea is to factorize the
given surjective morphism, after some transformation, through a
universal morphism associated to the minimal rational curves.
\end{abstract}


\medskip
2000MSC: 14J45, 32H02

\bigskip
\section{Introduction}

We will work over the complex numbers. A variety or (a manifold)
will be assumed to be irreducible except when we say `the variety of
minimal rational tangents', which may have finitely many components.
See Section 2 for the definition. For a complex manifold $Y$, $T(Y)$
denotes its tangent bundle and $T_y(Y)$ denotes the tangent space at
a point $y \in Y$.

\medskip
 For two projective varieties $X$
and $Y$, denote by
 ${\rm Hom}^s(Y, X)$  the space of surjective holomorphic maps $Y \to X$
 and by ${\rm Aut}_o(X)$ the identity component of the group of
 biregular automorphisms of $X$. In [HM3] and [HM4], the following result
 was proved.

 \medskip
 {\bf Theorem 1.1} {\it Let  $X$ be a Fano manifold of Picard number 1 whose variety
 of minimal rational tangents is non-linear or finite. Then for any projective
 variety $Y$,
each component of the reduction ${\rm Hom}^s(Y, X)_{red}$ is a
principal homogeneous space under the affine algebraic group ${\rm
Aut}_o(X).$}

\medskip
Theorem 1.1 was first proved for the rational homogeneous space
$X=G/P$ in [HM1]. It was proved when the variety of minimal rational
tangents has non-degenerate Gauss map in [HM2]. This was surpassed
by [HM4] which proves it when the variety of minimal rational
tangents is non-linear. The proofs in these three papers are of the
same nature. The proof when the variety of minimal rational tangents
is finite is quite different and  appeared in [HM3].

\medskip The condition that the variety of minimal rational tangents
is non-linear or finite holds for practically all known examples
of Fano manifolds of Picard number 1, except the projective space.
In fact, we have the following non-linearity conjecture:

\medskip
{\bf Conjecture 1.2} {\it Let $X$ be a Fano manifold of Picard
number 1 whose variety of minimal rational tangents is linear and
of positive dimension. Then $X$ is biregular to the projective
space.}

\medskip
There are some partial results toward Conjecture 1.2. For example,
it was proved for Fano manifolds of index $\geq \frac{\dim X
+3}{2}$ in [H1, Corollary 2.3].

\medskip
 For the projective
space, the assertion in Theorem 1.1 certainly does not hold. In
this sense, Theorem 1.1 is a reasonably satisfactory result,
except that it does not say whether ${\rm Hom}^s(Y, X)$ is
reduced. In other words, it does not address the unobstructedness
of infinitesimal deformations. The goal of this paper is precisely
to remedy this. Our main result is the following, which also gives
an alternative proof of Theorem 1.1.

\medskip
{\bf Theorem 1.3} {\it Let $X$ be a Fano manifold of Picard number 1
whose variety of minimal rational tangents is non-linear or finite
and let $Y$ be a projective variety. If $f: Y \rightarrow X$ is a
surjective morphism, then
$$H^0(Y, f^* T(X)) = f^* H^0(X, T(X)).$$
In particular, all deformations of surjective morphisms $Y \to X$
are unobstructed and each component of ${\rm Hom}^s(Y,X)$ is a
reduced principal homogeneous space of the affine algebraic group
${\rm Aut}_o(X)$. }

\medskip
Note that $H^0(Y, f^*T(X))$ is the Zariski tangent space to ${\rm
Hom}^s(Y, X)$ at $[f]$ and $H^0(X, T(X))$ is the Zariski tangent
space to ${\rm Aut}_o(X)$ at the identity. Thus the identity
$$H^0(Y, f^* T(X)) = f^* H^0(X, T(X))$$ implies that the natural
morphism ${\rm Aut}_o(X) \to {\rm Hom}^s(Y, X)$ sending each $g\in
{\rm Aut}_o(X) $ to $ g \circ f \in {\rm Hom}^s(Y, X)$ is bijective,
implying the last sentence of Theorem 1.3.

\medskip
The proof of Theorem 1.3 when the variety of minimal rational
tangents is non-linear is rather simple modulo the main result of
[HM4] on the tangent map. In retrospect, this proof is the
culmination of successive refinements of the arguments in [HM2] and
[HM4]. The final formulation is much simpler than the old proofs and
will be given in Section 2.

The difficult case is when the variety of minimal rational tangents
is finite. The
 key idea of the proof in that case is to show that, after
a certain transformation, the morphism $f: Y \to X$ can be
factorized through the universal morphism for the family of minimal
rational curves. This factorization is established in Section 4.
Combining this with an idea from [H1] on the behavior of minimal
rational curves near the branch locus of $f$  explained in Section
5, the proof is completed in Section 6 by using an argument in [H2].

\section{ Proof of Theorem 1.3 when the variety of minimal rational tangents of $X$ is non-linear}

\medskip Throughout this paper, we will denote by $X$
 a  Fano manifold of Picard number 1. We refer the readers to
[K] for basics on the space of rational curves on $X$. An
irreducible component ${\cal K}$ of the space of rational curves
on $X$  is called a {\it minimal component} if for a general point
$x \in X$, the subscheme ${\cal K}_x$ of ${\cal K}$ consisting of
members passing through $x$ is non-empty and complete. In this
case, the subvariety ${\cal C}_x$ of the projectivized tangent
space ${\P}T_x(X)$ consisting of the tangent directions at $x$ of
members of ${\cal K}_x$ is called the {\it variety of minimal
rational tangents} at $x$ (see [HM4] for more details). We say
that the variety of minimal rational tangents of $X$ is {\it
non-linear} if $\dim {\cal C}_x
>0$ and some component of ${\cal C}_x$ is  not a linear subspace
in $\P T_x(X)$. Otherwise, we say that the variety of minimal
rational tangents is linear.

For a general member $C$ of $\K$, the normalization $\nu: \P_1 \to C
\subset X$ is an immersion and $$\nu^*T(X) = \O(2) \oplus \O(1)^p
\oplus \O^q$$ where $p = \dim \C_x$ for a general $x \in X$ and $p+q
= \dim X -1$.  Denote by $H^0(C, T^*(X))$ the vector space
$$H^0(\P_1, \nu^*T^*(X)) = H^0(\P_1, \O^q).$$ For a non-singular point
$x \in C$, denote  by $H^0(C, T^*(X))_x \subset T_x^*(X)$ the
$q$-dimensional subspace of the cotangent space at $x$ given by
evaluating the elements of $H^0(C, T^*(X))$ at the point $x$.

\medskip
{\bf Proposition 2.1} {\it Let $X$ and $\K$ be as above and let $x
\in X$ be a general point. Let ${\cal K}_1$ be an irreducible
component of ${\cal K}_x$. Suppose that there exists a non-zero
vector $v \in T_x(X)$ annihilating $H^0(C, T^*(X))_x \subset
T_x^*(X)$ for any general member $C$ of ${\cal K}_1$. Then the
variety of minimal rational tangents of $X$ is linear.}

\medskip
{\it Proof}. From the irreducibility of $\K$, it suffices to show
that the component $\C_1$ of $\C_x$ corresponding to $\K_1$ is a
linear subspace. For a general member $C$ of ${\cal K}_1$, $x$ is a
non-singular point of $C$. Denote by
$$ \P H^0(C, T^*(X))_x \subset \P T_x^*(X)$$ the projectivization of
$H^0(C, T^*(X))_x$. The closure of the union of $\P H^0(C,
T^*(X))_x$ as $C$ varies over general points of ${\cal K}_1$, is the
dual variety of ${\cal C}_1 \subset \P T_x(X)$ by [HR, Corollary
2.2]. Thus the existence of $v$ implies that the dual variety of
${\cal C}_1$ is linearly degenerate in ${\bf P}T^*_x(X)$, i.e.,
${\cal C}_1$ is a cone. Thus Proposition 2.1 follows from [HM4,
Proposition 13], which says that ${\cal C}_1$ cannot be a cone
unless it is a linear subspace. $\Box$

\medskip
The next proposition is [HM2, Lemma 4.2].

\medskip
{\bf Proposition 2.2} {\it Let $X$ and $ \K$ be as above.  Let $Y$
be a projective variety and  $f: Y \rightarrow X$ be a generically
finite morphism of degree $>1$. Given a general member $C \subset X$
of ${\cal K}$, there exists a component $C'$ of $f^{-1}(C)$ such
that the restriction $f|_{C'} : C' \rightarrow C$ is finite of
degree $>1$.}

\medskip
{\bf Proposition 2.3} {\it In the situation of Proposition 2.2, let
$x \in C$ be a non-singular point and let $y_1, y_2 \in C'$ be two
distinct points with $f(y_1)=f(y_2) =x$. For a given $\sigma \in
H^0(Y, f^*T(X))$, regard its value $\sigma_{y_i} \in
(f^*T(X))_{y_i}$ as a vector in $T_x(X)$ for each $i=1,2$. Then
$\sigma_{y_1} - \sigma_{y_2} \in T_x(X)$ annihilates $H^0(C,
T^*(X))_x$. }

\medskip
{\it Proof}. Let $\nu: \hat{C} \to C$ be the normalization of $C$
and let $\varphi \in H^0(\hat{C}, \nu^*T^*(X))$ be  a section of the
 cotangent bundle of $X$ on $\hat{C}$. Let $\nu': \hat{C}' \to C'$ be the normalization
 of $C'$ and $\hat{f}: \hat{C}' \to \hat{C}$ be the lifting of $f$.
    Let $\varphi' \in H^0(\hat{C}',
(\nu \circ \hat{f})^*T^*(X))$ be the section induced by $\varphi$
and $\hat{\sigma} \in H^0( \hat{C}', ( f \circ \nu')^* T(X))$ be
the section induced by $\sigma$. Since $ \nu \circ \hat{f} = f
\circ \nu'$, the pairing $\varphi'(\hat{\sigma})$ is a holomorphic
function on $\hat{C}'$, hence is constant. It follows that
$\varphi'(\sigma_{y_1} ) = \varphi'(\sigma_{y_2}).$ Thus
$\sigma_{y_1} - \sigma_{y_2}$ is annihilates the evaluation of
$\varphi$ at $x$. $\Box$

\medskip
Now we can prove Theorem 1.3 when the variety of minimal rational
tangents of $X$ is non-linear.

\medskip
{\bf Proposition 2.4} {\it Let $X$ and $ \K$ be as above. Suppose
that there exists a surjective morphism $f:Y \to X$ from a
projective variety $Y$ with
$$H^0(Y, f^* T(X)) \neq f^* H^0(X, T(X)).$$ Then the variety of
minimal rational tangents of $X$ is linear.}

\medskip
{\it Proof}. Fix an element $\sigma \in H^0(Y, f^*T(X)) \setminus
f^* H^0(X, T(X)).$ For each $y \in Y$, let $\sigma_y \in
T_{f(y)}(X)$ be the corresponding tangent vector of $X$.
Associated to $\sigma$, we have the projective subvariety $\Sigma
\subset T(X)$ defined by
$$ \Sigma := \{ \sigma_y \in T_{f(y)} (X), y \in Y\}.$$
Since $\sigma \not\in f^*H^0(X, T(X))$, the natural projection $\pi:
\Sigma \to X$ is a finite morphism of degree $>1$ and $\sigma$
induces a natural section $\sigma'$ of $\pi^*T(X)$ with $\sigma'
\not\in \pi^* H^0(X, T(X))$. Thus replacing $Y$ by $\Sigma$, we may
assume that $f: Y \to X$ is a finite morphism and for any $x \in X$,
$$\sigma_{y_1} \neq \sigma_{y_2} \mbox{ as vectors in } T_x(X)
\mbox{ for each } y_1 \neq y_2 \in f^{-1}(x).$$

 Let $x$ be a general point of $X$ and $\K_1$ be an irreducible
component of $\K_x$. By Proposition 2.2, there exist two distinct
points $y_1, y_2 \in f^{-1}(x)$ such that for each general member
$C$ of $\K_1$, there exists an irreducible component $C'$ of
$f^{-1}(C)$ with $\{ y_1, y_2 \} \subset C'$. Then by Proposition
2.3, $H^0(C, T^*(X))_x$ is annihilated by $\sigma_{y_1}-
\sigma_{y_2}$ for all general members $C$ of $\K_1$. Applying
Proposition 2.1 with $v= \sigma_{y_1} - \sigma_{y_2}$, we conclude
that the variety of minimal rational tangents of $X$ is linear.
$\Box$

\section{Free curves with trivial normal bundle}

It is convenient to introduce the following notion. Let $Y$ be a
projective manifold of dimension $n$ and $C \subset Y$ be an
irreducible curve. We say that $C$ is a {\it free curve with
trivial normal bundle} if the following holds.

(i) Under the normalization $\nu: \hat{C} \to C$, we have an exact
sequence of vector bundles on $\hat{C}$
$$0 \longrightarrow T(\hat{C}) \longrightarrow \nu^* T(Y)
\longrightarrow N_C \longrightarrow 0$$  where the second arrow is
the differential of $\nu: \hat{C} \to Y$ and  $N_C$ is a trivial
bundle of rank $=(n-1)$ on $\hat{C}$.

(ii) Deformations of $C$ with constant geometric genus cover an
open subset of $X$.

\medskip The germ of the space of deformations of $C$ with constant
geometric genus must have dimension $\geq n-1$. The Zariski
tangent space to this space at the point corresponding to $C$ is
$H^0(\hat{C}, N_C)$, which has dimension $n-1$ from the triviality
of the normal bundle. Thus the germ of this space of deformations
of $C$, which we denote by ${\cal M}_C$, is non-singular. The
following is obvious from the deformation theory of submanifolds.

\medskip
{\bf Proposition 3.1} {\it Let $ C \subset Y$ be a free curve with
trivial normal bundle. Let $\vartheta$ be a germ of
nowhere-vanishing holomorphic vector fields on ${\cal M}_C$ given by
some element
$$\vartheta_{[C_s]} \in H^0(\hat{C}_s, N_{C_s}) \mbox{ for each }
[C_s] \in {\cal M}_C.$$ Denote by $\Delta$ the complex unit disc.
The integral curve of $\vartheta$ through $[C]$ defines a
deformation $\{ [C_t] \in {\cal M}_C, t \in \Delta, C= C_0\}$ of
$C$. Let $x \in C$ be a non-singular point. Suppose there exists a
germ $\theta$ of holomorphic vector fields of $Y$ at $x$ such that
$\theta$ modulo $T(C_s) $ agrees with the germ of
$\vartheta_{C_s}$ at $x$ for each $[C_s] \in {\cal M}_C$. Then the
integral curve of $\theta$ through $x$ defines a deformation $\{
x_t \in Y, t \in \Delta, x= x_0\} $ of $x$ such that $x_t \in C_t$
for each $t$, up to reparametrization.}

\medskip
From now throughout the rest of this paper, we will fix a Fano
manifold $X$ of Picard number 1 and a minimal component $\K$ such
that the variety of minimal rational tangents at a general point is
finite. Then  a general member $C$ of $\K$ is a free curve with
trivial normal bundle and the germ ${\cal M}_C$ can be realized by
an open neighborhood of $[C] \in \K$. By desingularizing the
universal family  over $\K$ (see [K, II.2.12] for the definition of
the universal family), we have the following. The proof, which is
quite standard, will be omitted.

\medskip
{\bf Proposition 3.2} {\it There exist a  projective manifold $X'$
with a generically finite morphism $\mu:X' \rightarrow X$ of
degree $>1$ and a proper surjective morphism $\rho: X' \rightarrow
Z$ onto a projective manifold $Z$ with the following properties.

(a)  $\rho$ is a $\P_1$-bundle over a Zariski open dense subset
$Z_o\subset Z$.

(b)  $\mu$ is unramified on $\rho^{-1}(Z_o)$.

(c)  Each member of $\K_x$ for a general $x \in X$ is the image of
 a fiber of $\rho$ through $\mu^{-1}(x)$.

(d) For each $\zeta \in Z_o$,  $\mu|_{\rho^{-1}(\zeta)}$ is the
normalization of $P_{\zeta} := \mu(\rho^{-1}(\zeta)).$

 (e) For two distinct points $\zeta_1 \neq \zeta_2 \in Z_o$, the two
curves $P_{\zeta_1}$ and $P_{\zeta_2}$ are distinct.}

\medskip
Let us denote by $T^{\rho} \subset T(\rho^{-1}(Z_o))$ the relative
tangent bundle of $\rho$ over $\rho^{-1}(Z_o)$.  Let ${\cal C}
\subset {\bf P}T(X)$ be the closure of the union of ${\cal C}_x$'s
for general points $x \in X$. Let $\hat{\C} \subset T(X)$ be the
cone over $\C$. Denote by $O \subset T(X)$ the zero section and by
$\pi: T(X) \to X$ the natural projection. The following is
immediate.

\medskip
{\bf Proposition 3.3} {\it In the setting of Proposition 3.2, there
exists a Zariski open dense subset $X_o \subset X$ such that
$\mu^{-1}(X_o) \subset \rho^{-1}(Z_o)$ and the restriction of $\pi$
to $(\hat{\C} \setminus O) \cap \pi^{-1}(X_o)$ is a smooth morphism.
  For each point $x \in X_o$ and $\mu^{-1}(x) = \{
x_1, \ldots, x_m\}, m= $ degree of $\mu$, we have a disjoint union
$$ \pi^{-1}(x) \cap (\hat{\C} \setminus O) = d \mu (T^{\rho}_{x_1} \setminus \{0\}) \cup \cdots \cup d \mu
(T^{\rho}_{x_m} \setminus \{ 0 \}). $$ In particular, we have a
natural smooth morphism $$ \chi: \; (\hat{\C} \setminus O) \cap
\pi^{-1}(X_o) \longrightarrow \mu^{-1}(X_o)$$ such that $\pi = \mu
\circ \chi$ on $(\hat{\C} \setminus O) \cap \pi^{-1}(X_o)$.}

\medskip
{\bf Proposition 3.4} {\it Let $Y$ be a projective manifold and
$f: Y \to X$ be a generically finite surjective morphism. For a
general member $C \subset X$ of $\K$, $C$ intersects each
component of the branch divisor of $f$ transversally. Each
irreducible component $C'$ of $f^{-1}(C)$ is a free curve with
trivial normal bundle and when $\hat{C}$ (resp. $\hat{C}'$) is the
normalization of $C$ (resp. $C'$) and $\hat{f}: \hat{C}' \to
\hat{C}$ is the morphism induced by $f$, there are canonical
isomorphisms
$$T_{[C']}({\cal M}_{C'}):= H^0( \hat{C}', N_{C'}) \cong H^0(\hat{C}', \hat{f}^*N_C))
\cong H^0(\hat{C}, N_C) =: T_{[C]}({\cal M}_C)$$ and a biholomorphic
equivalence of germs ${\cal M}_C \cong {\cal M}_{C'}$. }

\medskip
{\it Proof}. That $C$ intersects the branch divisor transversally
is obvious from Proposition 3.2 (b).  The fact that $C'$ is a free
curve with trivial normal bundle is precisely [HM3, Proposition
6]. The canonical isomorphisms and the equivalence of germs are
obvious from the isomorphism of two trivial vector bundles $
N_{C'}  \cong \hat{f}^* N_C$ induced by the differential $df: T(Y)
\to \hat{f}^*T(X)$. $\Box$

\medskip
{\bf Proposition 3.5} {\it Let $Y$ be a projective variety and $f:
Y \to X$ be a generically finite surjective morphism. Let $C$ be a
general member of $\K$ and let $C'$ be a component of $f^{-1}(C)$.
Pick a non-singular point $x \in C$ outside the branch loci. Let
$\sigma \in H^0(Y, f^*T(X))$. For any two points $y_1, y_2 \in
f^{-1}(x) \cap C'$, regard $\sigma_{y_1}$ and $ \sigma_{y_2}$ as
vectors in $T_x(X)$. Then $\sigma_{y_1} - \sigma_{y_2}  \in
T_x(C)$. In particular,   $\sigma$ induces a unique element in
$H^0(\hat{C}, N_{C})$, up to a choice of $C'$.}

\medskip
{\it Proof}. This is a consequence of Proposition 2.3. Since $C$ has
 trivial normal bundle,  $H^0(C, T^*(X))_x$ is the  conormal space
of $C$ at $x$. Thus $\sigma_{y_1} - \sigma_{y_2} \in T_x(C)$. $\Box$

\section{Factorization through $\mu$}

In the setting of Theorem 1.3, given a section $\sigma \in H^0(Y,
f^*T(X))$, the values of $\sigma$ define a projective variety in
$T(X)$ dominant over $X$, as explained in the proof of Proposition
2.4. In fact, Theorem 1.3 is equivalent to the statement that a
projective variety in $T(X)$ dominant over $X$ must be a section of
$T(X)$. In other words, we have to prove that there do not exist
projective varieties of $T(X)$ which have degree $>1$ over $X$. The
goal of this section is to show that given a projective variety
$\Sigma \subset T(X)$ of degree $>1$ over $X$, the difference
transform of $\Sigma$ contains an irreducible
 component  that has very special properties with respect to the morphisms
$\mu, \rho$ of Proposition 3.2. It should be mentioned that all
the propositions proved from now on, except Proposition 5.1, are
under the assumption of the existence of $\Sigma $ of degree $>1$,
which will lead eventually to contradiction. In this sense all
these propositions are of hypothetical nature.

\medskip
{\bf Proposition 4.1} {\it  Suppose there exists a projective
variety $\Sigma \subset T(X)$ which is dominant over $X$ of degree
$>1$. Let  $T(X) \times_X T(X)$ be the fiber product of two copies
of the projection $\pi: T(X) \to X$ and let
$$\Sigma \times_X \Sigma \subset T(X) \times_X T(X)$$ be the fiber
product of two copies of $\pi|_{\Sigma}: \Sigma \to X$. Then there
exists an irreducible component $\Sigma^{\sharp}$ of $\Sigma
\times_X \Sigma$ with the following property:
 for a general $\zeta \in Z_o$ and a general point $x \in
 P_{\zeta}$, there exists an irreducible component $P'_{\zeta}$ of $\pi^{-1}(P_{\zeta}) \cap
\Sigma$ and two distinct points $a_1 \neq a_2 \in P'_{\zeta} \cap
\pi^{-1}(x)$ such that  $\Sigma^{\sharp}$, regarded as a subvariety
of
  $\Sigma \times \Sigma$, contains the point $(a_1, a_2)$. }

\medskip
{\it Proof}. For a general $\zeta \in Z_o$, there exists an
irreducible component $P'_{\zeta}$ of $\pi^{-1}(P_{\zeta}) \cap
\Sigma$ such that the projection $P'_{\zeta} \to P_{\zeta}$ is
finite of degree $>1$ by Proposition 2.2. Thus we can choose two
$a_1 \neq a_2$ on $P'_{\zeta}$ over $x \in P_{\zeta}$. The point
$(a_1, a_2) \in \Sigma \times \Sigma$ lies in $\Sigma \times_X
\Sigma.$ From the generality of the choice of $\zeta$ and $x$, there
is a unique component $\Sigma^{\sharp}$ containing $(a_1, a_2)$.
Certainly, $\Sigma^{\sharp}$ satisfies the required property from
the irreducibility of $Z_o$. $\Box$

\medskip
{\bf Proposition 4.2} {\it In the situation of Proposition 4.1, let
$\delta: T(X) \times_X T(X) \to T(X)$ be the difference morphism
defined by
$$\delta( v_1, v_2) := v_1 - v_2 \; \mbox{ for } \; v_1, v_2 \in T_x(X) \; \mbox{ for } x \in X.$$
Then in the notation of Proposition 3.3, $$\delta (\Sigma^{\sharp})
\subset \hat{\C},\;\;  \delta(\Sigma^{\sharp}) \not\subset O,$$ and
the dominant rational map $\chi^{\sharp}: \delta(\Sigma^{\sharp})
\longrightarrow X',$ induced by the morphism $\chi$, is generically
finite. }

\medskip
{\it Proof}. We will apply Proposition 3.5 with  $Y = \Sigma, f=
\pi|_{\Sigma}$ and $C = P_{\zeta}$. There is a tautological
section $\sigma \in H^0(Y, f^*T(X))$ defined by
$$ \sigma_a = a \in T_x(X) \; \mbox{ for each } a \in \Sigma \cap
T_x(X).$$ By Proposition 3.5, $$a_1 - a_2 \; \in  \;
T_x(P_{\zeta}) \; \subset \; \hat{\C}.$$ As $\zeta$ varies over
general points of $Z_o$, the element $a_1-a_2$ varies over an open
subset in the irreducible $\delta( \Sigma^{\sharp})$. It follows
that $\delta(\Sigma^{\sharp}) \subset \hat{\C}$. Since $a_1 \neq
a_2$, $\delta(\Sigma^{\sharp})$ is not contained in the zero
section $O$. The dominant rational map $\chi^{\sharp}$ is
certainly generically finite. $\Box$

\medskip
{\bf Proposition 4.3} {\it  In the situation of Proposition 4.2,
there
  exists a  projective manifold $\tilde{\Sigma}$, a generically finite morphism $g:
\tilde{\Sigma} \to X'$ and a section $\theta \in
H^{0}(\tilde{\Sigma}, (\mu \circ g)^* T(X))$ with the following
properties.

(1) For a general point $x \in X$ and any two distinct $y_1, y_2 \in
(\mu \circ g)^{-1}(x)$, $\theta_{y_1} \neq \theta_{y_2}$ as vectors
in $T_x(X)$,

(2) For a general point $x \in X$ and any $y \in (\mu \circ
g)^{-1}(x)$, $\theta_y$ regarded as   a vector in $T_{g(y)}(X') =
T_x(X)$, belongs to $ T^{\rho}_{g(y)}$  where $T^{\rho}$ is as in
Proposition 3.3.}

\medskip
{\it Proof}. Choose a desingularization $\alpha: \tilde{\Sigma}
\to \delta(\Sigma^{\sharp})$ which eliminates the indeterminacy of
the generically finite rational map $\chi^{\sharp}$  such that $
\chi^{\sharp} \circ \alpha$ defines a generically finite morphism
$g: \tilde{\Sigma} \to X'$.  Denote by $\tau$ the natural
projection $\delta(\Sigma^{\sharp}) \to X$. Then $\tau \circ
\alpha = \mu \circ g$.  Since $\delta(\Sigma^{\sharp}) \subset
T(X)$, there exists a tautological section $\kappa \in
H^0(\delta(\Sigma^{\sharp}), \tau^* T(X))$ defined by $\kappa(a) =
a \in T_{\tau(a)}(X)$ for each $a \in \delta(\Sigma^{\sharp})$.
Let
$$\theta \in H^0( \tilde{\Sigma}, (\mu \circ g)^* T(X)) = H^0(
\tilde{\Sigma}, (\tau \circ \alpha)^* T(X))
$$ be the pull-back of $\kappa$ by $\alpha$. Then $\theta$ satisfies
 property (1), because $\alpha$ is birational and the
tautological section $\kappa$ satisfies an analog of (1). It
satisfies  property (2) from $\delta (\Sigma^{\sharp}) \subset
\hat{\C}.$ $\Box$

\section{Univalence of $\K$ on the branch divisor of $\mu$}

In the setting of Proposition 3.2, we say that $\K$ is {\it
univalent on an irreducible hypersurface} $B \subset X$ if  (i)
there exists only one irreducible component $E$ of $\mu^{-1}(B)$
that is dominant over both $Z$ and $B$, and (ii) the morphism
$\mu|_E: E \to B$ is birational.  This is equivalent to saying
that at a general point $z \in B$, there exists exactly one member
$C$ of $\K$ passing through $z$ with $C \not\subset B$ and $C$ is
non-singular at $z$.

\medskip
The following is  essentially the same as  [H1, Proposition 3.2].

\medskip
{\bf Proposition 5.1} {\it  In the setting of Proposition 3.2,
suppose that there exists an irreducible hypersurface $B \subset
X$, such that $\K$ is not univalent on $B$. Then given  a general
point $x \in B$ and an open neighborhood $W\subset X$ of $x$,
there exists a point $y \in W$ and two distinct points $\zeta_1,
\zeta_2 \in Z_o$ with $y \in P_{\zeta_1} \cap P_{\zeta_2}$ and
$T_y(P_{\zeta_1}) \neq T_y(P_{\zeta_2})$ such that irreducible
components of $W \cap P_{\zeta_1}$ and $W \cap P_{\zeta_2}$
through $y$ intersect $B$ transversally at some point of $B \cap
W$. }

\medskip
{\it Proof}.  By assumption, there exist a union $D$ of components
of $\mu^{-1}(B)$ each of which is dominant over $Z$ and $B$, and
the morphism $\mu|_D: D  \to B$ has degree $>1$.  Let $y_1, y_2$
be two distinct points of $\mu^{-1}(x) \cap D$. Since $x$ is
general, both $\rho(y_1)$ and $\rho(y_2)$ lie in $Z_o$. There
exist open neighborhoods $W_1 \subset \rho^{-1}(Z_o)$ of $y_1$,
$W_2 \subset \rho^{-1}(Z_o)$ of $y_2$ and $W_0 \subset W$ of $x$
with the following properties

(1) $\mu(W_1) = \mu(W_2) = W_0$,

(2)  $\mu|_{W_1}$ and $\mu|_{W_2}$ are biholomorphic,

(3)  $W_1 \cap D$ and $W_2 \cap D$ are non-singular and
transversal to the fibers of $\rho$.

There exist open neighborhood $W_1' \subset W_1$ of $y_1$ and $W'_2
\subset W_2$ of $y_2$ such that for any $y \in W_1'$ (resp. $y \in
W_2'$) $\rho^{-1}(\rho(y)) \cap W_1'$ (resp. $\rho^{-1}(\rho(y))
\cap W_2'$) is connected. Let $y$  be a general point in $ \mu(W_1')
\cap \mu(W_2')$. Let $y_1' = W_1' \cap \mu^{-1}(y)$ and $y_2' = W_2'
\cap \mu^{-1}(y).$  Then $\zeta_1 := \rho(y_1')$ and $\zeta_2 :=
\rho(y_2')$ give the desired two distinct points. $\Box$

\medskip
The idea of the proof of the following proposition is  the same as
that of [H1, Proposition 3.3].

\medskip
{\bf Proposition 5.2} {\it In the setting of Proposition 4.3, $\K$
is univalent on each component $B$ of the branch divisor of
$\mu$.}

\medskip
{\it Proof}. Suppose that $\K$ is not univalent on some component
$B$.  Set $Y:= \tilde{\Sigma}$ and $f= \mu \circ g$. Then $B$ is a
component of the branch divisor of $f: Y \to X$. Let $R \subset Y$
be an irreducible component of the ramification divisor of $f$
such that $B =f(R)$. Let $z \in R$ be a general point and let $r$
be the local sheeting number of $f$ at $z$. We can choose a
holomorphic coordinate neighborhood $V$ of $z$ with coordinates
$(w_1, \ldots, w_n)$ at $z$ and a holomorphic coordinate
neighborhood $W$ of $f(z)$ with coordinates $(z_1, \ldots, z_n)$
such that $B \cap W$ is defined by $z_n=0$ and  $f$ is given by
$$z_1 = w_1, \ldots, z_{n-1} = w_{n-1}, z_n = w_n^r.$$
Let $x \in W\setminus B$ and  $\zeta_1, \zeta_2 \in Z_o$ be as in
Proposition 5.1. Setting $C_1=P_{\zeta_1}$ (resp. $C_2
=P_{\zeta_2}$), an easy coordinate computation in the above
coordinate systems (see e.g. [HM3, p.636, Lemma1]) shows that there
exists a unique irreducible component $C'_1$ (resp. $C'_2$) of
$f^{-1}(C_1)$ (resp. $f^{-1}(C_2)$) intersecting $V$ such that an
irreducible component of $C'_1 \cap V$ (resp. $C'_2 \cap V$)
contains $f^{-1}(x) \cap V$. In particular, $C'_1 \cap C'_2$
contains the $r$ distinct points $f^{-1}(x) \cap V$. Let $y_1 \neq
y_2$ be two distinct points in $f^{-1}(x) \cap V \cap C'_1 \cap
C'_2$. Applying Proposition 3.5 to $C'_1$ and $C_1$,
$$ \theta_{y_1} - \theta_{y_2} \in T_{x_o}(C_1).$$ Applying Proposition
3.5 to $C'_2$ and $C_2$, $$ \theta_{y_1} - \theta_{y_2} \in
T_{x_o}(C_2).$$ Since $T_{x_o}(C_1) \cap T_{x_o}(C_2) =0$, we get
$\theta_{y_1} = \theta_{y_2}$, a contradiction to Proposition 4.3
(1). $\Box$

\medskip
{\bf Proposition 5.3} {\it In the setting of Proposition 5.2, let
$B \subset X$ be a component of the branch divisor of $\mu$ and
let $D$ be the unique irreducible component of $\mu^{-1}(B)$ which
is dominant over $Z$  and $B$. Then for a general member $C$ of
$\K$, any component $C'$ of $\mu^{-1}(C)$ which is finite over $C$
of degree $>1$ is disjoint from $D$.}

\medskip {\it Proof}.
Suppose not. Since $C'$ is a free curve with trivial normal bundle
by Proposition 3.4, we may assume that $C'$ intersects $D$ at a
general point $x'$ of $D$. Then through a general point $x'$ of
$D$,
  we have two distinct curves, $C'$ and a fiber of
$\rho$, neither of which are contained in $D$. Since $\mu$ is
unramified at $x'$ by Proposition 3.2 (b), the images of these
curves under $\mu$ are of the form $P_{\zeta_1}, P_{\zeta_2}$ with
$\zeta_1 \neq \zeta_2$. Since these two curves pass through $x =
\mu(x')$, which  is a general point of $B$, ${\cal K}$ is not
univalent on $B$, a contradiction to Proposition 5.2. $\Box$

\section{Completion of the proof of Theorem 1.3}

In this section, we will complete the proof of Theorem 1.3.  The
strategy is to establish some analogs of [H2, Section 5].

\medskip
{\bf Proposition 6.1} {\it In the setting of Proposition 4.3, let
$C \subset X$ be a general member of $\K$. By Proposition 2.2,
there exists a component $C' $  of $\mu^{-1}(C)$ such that
$\mu|_{C'}: C' \to C$ is finite of degree $>1$, and $C'$ is a free
curve with trivial normal bundle by Proposition 3.4.  Fix  a
choice of a component $C^{\flat}$ of $g^{-1}(C')$. Then $\theta$
modulo $T(\hat{C}')$ defines  a non-zero element $\vartheta_{C'}
\in H^0(\hat{C}', N_{C'})$. }

\medskip
{\it Proof}. By Proposition 3.5, $C^{\flat}$ determines a unique
element in $H^0(\hat{C}, N_C)$. By the isomorphism in Proposition
3.4, this determines an element $\vartheta_{C'} \in H^0(\hat{C}',
N_{C'}).$ It cannot be zero because of Proposition 4.3 (2) and the
generality of $C$. $\Box$

\medskip
The next proposition is an analog of [H2, Lemma 5.5], although their
proofs are of different nature.

\medskip
 {\bf Proposition 6.2} {\it In the setting of Proposition
6.1, denoting by $\Delta$ the complex unit disc, there exist a
family of members of ${\cal K}$
$$\{C_t, t \in \Delta, C=C_0 \} $$ and the associated deformation
$$\{ C'_t, t \in \Delta, C'=C'_0\}$$  such that for each $t \in
\Delta$,

(i) $C'_t$ is a component of $\mu^{-1}(C_t)$;

(ii)  $\mu|_{C'_t}: C'_t \rightarrow C_t$ is finite of degree $>1$;

(iii) $\rho(C'_t) = \rho(C')$.}

\medskip
{\it Proof}. For a deformation $[C'_s] \in {\cal M}_{C'}$ of $C'$,
we get a deformation $[C^{\flat}_s] \in {\cal M}_{C^{\flat}}$ with
$C^{\flat}_s \subset g^{-1}(C'_s).$ By applying Proposition 6.1 to
$C^{\flat}_s$, we get an element $\vartheta_{C'_s} \in H^0(C'_s,
N_{C'_s}).$ Thus the choice of $C^{\flat}$ determines a germ of
holomorphic vector fields $\vartheta$ on ${\cal M}_{C'}$. By
Proposition 4.3 (2), this is a germ of non-vanishing vector fields.
Let
$$\{ C'_t, t \in \Delta, C'=C'_0\} = \{ [C'_t] \in {\cal M}_{C'}, t \in \Delta \}$$ be a local analytic arc
integrating the vector field $\vartheta$. (i) and (ii) are obvious
from the local equivalence ${\cal M}_C \cong {\cal M}_{C'}$ in
Proposition 3.4. It suffices to check (iii). Let $x \in C'$ be a
general point and $y \in C^{\flat}$ be the point over $x$. Then
the germ of holomorphic vector fields defined by $\theta$ at $y$
induces a germ $\theta'$ of holomorphic vector fields in a
neighborhood of $x$. Applying 3.1, we see that the integral curve
of $\theta'$ through $x$, $\{x_t \in X', t \in \Delta\}$ with $x =
x_0$, satisfies $x_t \in C'_t$ up to reparametrization.  Since
$\theta'$ is a section of $T^{\rho}$ by Proposition 4.3 (2), $x_t
\in \rho^{-1}(\rho(x))$. It follows that
$$  \;\; \rho^{-1}(z) \cap C'_t \neq \emptyset \mbox{ for general
} z \in \rho(C') \mbox{ and each } t \in \Delta.$$ This implies
(iii). $\Box$

\medskip   The proof of the next
proposition is, modulo Proposition 5.3 and Proposition 6.2,  almost
verbatim that of [H2, Proposition 5.3]. Since the terms and the
notation are slightly different, we reproduce the proof for the
reader's convenience.

\medskip  {\bf Proposition 6.3} {\it Let us assume the situation of
Proposition 6.1.  Given $C$ and $C'$ as in Proposition 6.1,  let
$$h: \hat{C}' \longrightarrow \widehat{\rho(C')}$$ be the lift of
$$\rho|_{C'}: C' \longrightarrow \rho(C')$$ to the normalizations
of $C'$ and $\rho(C')$. Then $h$ has a ramification point $z \in
\hat{C}'$ such that the image of $h(z)$ in $\rho(C')$ lies in
$Z_o$.}

\medskip
{\it Proof}.
 Suppose not. Then $h$ is unramified over $\rho(C') \cap Z_o$.
  Let us use the deformation $C_t$ constructed in Proposition 6.2.
   By the generality of $C$, we may assume that for each $t \in
\Delta$ the holomorphic map $$h_t: \hat{C}'_t \rightarrow
\widehat{\rho(C'_t)} = \widehat{\rho(C')}, \;\;  h_0 =h,$$ which
is the lift of $\rho|_{C'_t}$ to the  normalization, is unramified
over $\rho(C') \cap Z_o$. Since $h_t$ is a continuous family of
coverings of the Riemann surface $\widehat{\rho(C')}$ with fixed
branch locus, we can find a biholomorphic map
$$ (\clubsuit) \;\;\; \psi_t: \hat{C}' \rightarrow \hat{C}'_t, \;\; \psi_0 = {\rm Id}_{\hat{C}'} \mbox{ with } h = h_t \circ
\psi_t,$$ which depends holomorphically on $t$ (e.g. [S, p. 32,
Corollary 1].).

By Proposition 6.2 (ii), there are at least two distinct points in
$\hat{C}_t$, say $a_t \neq b_t \in \hat{C}_t$, such that the
corresponding points in $C_t$ lie in the branch divisor of $\mu$
in $X$. Let $\{ 0, \infty\} \subset \P_1$ be two distinct points
on the projective line. We can choose a family of biholomorphic
maps $\{ \sigma_t: \hat{C}_t \rightarrow \P_1, \; t \in \Delta \}$
such that $\sigma_t(a_t) = 0$ and $ \sigma_t(b_t) = \infty$ for
each $t \in \Delta$.  Denote by
 $\mu_t: \hat{C}'_t \rightarrow \hat{C}_t$  the lift of
$\mu|_{C'_t}$ to the normalization. Then
 $$ \{ \varphi_t: \hat{C}' \longrightarrow
\P_1, \;\;\; \varphi_t := \sigma_t \circ \mu_t \circ \psi_t, \; t
\in \Delta \}$$ is a family of meromorphic functions on the compact
Riemann surface $\hat{C}'$.

By Proposition 5.3, for each component $E$ of the branch divisor of
$\mu$, the intersection of $C'_t$ with $\mu^{-1}(E)$ has a fixed
image in $\rho(C')=\rho(C'_t)$, independent of $t \in \Delta$. This
implies that there is a finite subset $Q \subset
\widehat{\rho(C')}$, independent of $t$, such that
$$\mu_t^{-1}(a_t) \cup \mu_t^{-1}(b_t) \subset h_t^{-1}(Q)$$ for any $t
\in \Delta$. Then $$ \varphi_t^{-1}(0) = \psi^{-1}_t \circ
\mu_t^{-1} \circ \sigma_t^{-1}(0) = \psi^{-1}_t(\mu_t^{-1}(a_t))
\subset  \psi_t^{-1}(h_t^{-1}(Q)) $$ for all $t \in \Delta$. Since
$\psi_t^{-1}(h_t^{-1}(Q))= h^{-1}(Q)$ by the choice of $\psi_t$ in
$(\clubsuit)$, $ \varphi_t^{-1}(0) \subset h^{-1}(Q)$ for any $t
\in \Delta$. Consequently,  $\varphi_t^{-1}(0) =
\varphi_0^{-1}(0)$ for all $t \in \Delta$. By the same argument we
get $\varphi_t^{-1}(\infty) = \varphi_0^{-1}(\infty)$ for all $t
\in \Delta$. In other words, the family of meromorphic functions
$\varphi_t$ have the same zeroes and the same poles on the Riemann
surface $\hat{C}'$.  It follows that for any $w_1, w_2 \in
\hat{C}'$ and any $t \in \Delta$, $$ (\diamondsuit) \;\;\;\;
\varphi_t(w_1) = \varphi_t(w_2) \mbox{ if and only if }
\varphi_0(w_1) = \varphi_0(w_2).$$

Since $\mu|_{C'}$ is  finite of degree $>1$ by our assumption, we
can choose two points $\alpha \neq \beta \in \hat{C}'$ such that
$\varphi_0(\alpha) = \varphi_0(\beta).$  Furthermore, denoting by
$\bar{\alpha}\in \rho(C')$ (resp. $\bar{\beta} \in \rho(C')$)  the
point corresponding to $h_0(\alpha) \in \widehat{\rho(C')}$ (resp.
$h_0(\beta)\in \widehat{\rho(C')}$) under the normalization, we may
assume that $$ (\heartsuit) \;\;\; \bar{\alpha} \mbox{ and }
\bar{\beta} \mbox{ are two distinct points in } Z_o.$$ From
$(\diamondsuit)$, we have  $\varphi_t(\alpha) = \varphi_t(\beta)$
for all $t \in \Delta$. Since $\varphi_t = \sigma_t \circ \mu_t
\circ \psi_t$ and $\sigma_t$ is biholomorphic, we see that $$
(\spadesuit) \;\;\;\;  \mu_t \circ \psi_t(\alpha) = \mu_t \circ
\psi_t(\beta) \mbox{ for all } t \in \Delta.$$

 Denote by $$ \alpha_t \in C'_t \subset X' \;\;\; \mbox{ (resp. }
 \beta_t \in C'_t \subset X') $$ the point corresponding to
  $\psi_t(\alpha) \in \hat{C}'_t$ (resp.
 $\psi_t(\beta) \in \hat{C}'_t$)  under the normalization. Then the
 locus $$A := \{ \alpha_t \in X', t \in \Delta \} \;\;\; \mbox{ (resp.
  } B := \{ \beta_t \in X', t \in \Delta \})$$
covers a non-empty  open subset in the fibre
$\rho^{-1}(\bar{\alpha})$ (resp. $\rho^{-1}(\bar{\beta})$). Thus
$\mu(A)$ (resp. $\mu(B)$) covers a non-empty open subset in $$
P_{\bar{\alpha}} := \mu(\rho^{-1}(\bar{\alpha})) \;\;\; \mbox{
(resp. } P_{\bar{\beta}} :=\mu( \rho^{-1}(\bar{\beta} ))).$$ Since
$\mu(A)$ (resp. $\mu(B)$) is the locus of points corresponding to
$\mu_t \circ \psi_t (\alpha)$ (resp. $\mu_t \circ \psi_t (\beta)$)
by the normalization $\hat{C}_t \rightarrow C_t$, the equality
$(\spadesuit)$ above implies that $\mu(A) = \mu(B)$. Consequently,
$$P_{\bar{\alpha}}= P_{\bar{\beta}},$$ a contradiction to  Proposition 3.2
via $(\heartsuit)$. $\Box$

\medskip
Now we are ready to finish the proof of Theorem 1.3 as follows.

\medskip
{\it End of the proof of Theorem 1.3}. As explained at the beginning
of Section 4, we may assume the situation of Proposition 4.3 and get
a contradiction.  From Proposition 6.3, let $z \in C'$ be the image
of a ramification point of $h$ such that $\rho(z) \in Z_o$. Then a
component of the germ of $C'$ at $z$ must be tangent to
$T^{\rho}_z$.  Choose $C^{\flat}$ as in Proposition 6.1. The value
of $\theta$ at a point of $C^{\flat}$ over $z$ determines $\theta_z
\in T_z(X')$ which is in $T^{\rho}_z$ by Proposition 4.3 (2). Thus
$\theta_z$ is tangent to a component of the germ of $C'$ at $z$.
This means that the non-zero element $\vartheta_{C'} \in
H^0(\hat{C}', N_{C'})$ in Proposition 6.1 vanishes at $z$, a
contradiction to the triviality of $N_{C'}$. $\Box$

\bigskip
{\bf Acknowledgment} An essential idea for this work was obtained
during my visit to Fudan University in December, 2007. I would
like to thank Ngaiming Mok and Yuxin Dong for the invitation and
the hospitality.

\bigskip
\bigskip
{\bf References}
\medskip

 [H1] Hwang, J.-M.: Deformation of
holomorphic maps onto Fano manifolds of second and fourth Betti
numbers 1. Ann. Inst. Fourier {\bf 57} (2007) 815-823

[H2] Hwang, J.-M.: Base manifolds for fibrations of projective
irreducible symplectic manifolds. Invent. math. {\bf 174} (2008)
625-644

[HM1]  Hwang, J.-M. and Mok, N.: Varieties of minimal rational
tangents on uniruled manifolds. in {\it Several Complex Variables},
ed. by M. Schneider and Y.-T. Siu, MSRI Publications 37, Cambridge
University Press (2000) 351-389

[HM2]  Hwang, J.-M. and Mok, N.:  Cartan-Fubini type extension of
holomorphic maps for Fano manifolds of Picard number 1. Journal
Math. Pures Appl. {\bf 80}  (2001) 563-575

[HM3] Hwang, J.-M. and Mok, N.: Finite morphisms onto Fano manifolds
of Picard number 1 which have rational curves with trivial normal
bundles. J. Alg. Geom. {\bf 12} (2003) 627-651

 [HM4] Hwang, J.-M. and Mok,
N.: Birationality of the tangent map for minimal rational curves.
Asian J. Math. {\bf 8}, {\it Special issue dedicated to Yum-Tong
Siu}, (2004) 51-64

 [HR] Hwang, J.-M. and Ramanan,
S.: Hecke curves and Hitchin discriminant. Ann. scient. Ec. Norm.
Sup. {\bf 37} (2004) 801-817

[K] Koll\'ar, J.: Rational curves on algebraic varieties. Erg. d.
Math. 3 Folge, Band 32. Springer Verlag (1996)

[S] Shokurov, V. V.: Riemann surfaces and algebraic curves. In
{\it Algebraic curves, algebraic manifolds and schemes}, Springer
Verlag (1998)

\vspace{10mm}

Jun-Muk Hwang

Korea Institute for Advanced Study

207-43 Cheongnyangni-dong

Seoul 130-722, Korea

 jmhwang@kias.re.kr

\end{document}